\documentclass[12pt,reqno]{amsart}
\usepackage{amsmath,amssymb,amsfonts,amscd,latexsym,amsthm,mathrsfs}
\usepackage[usenames]{color}
\usepackage[unicode]{hyperref}
\newtheorem{theorem}{Theorem}

\theoremstyle{remark}
\newtheorem*{acknowledgements}{\bf Acknowledgements}

\let\eps\varepsilon
\renewcommand{\d}{{\mathrm d}}

\begin{document}

\title[Linear independence for Chowla--Selberg periods]{Linear independence measures \\ for Chowla--Selberg periods}

\author{Wadim Zudilin}
\address{Department of Mathematics, IMAPP, Radboud University, PO Box 9010, 6500~GL Nij\-me\-gen, Netherlands}
\urladdr{https://www.math.ru.nl/~wzudilin/}

\date{8 October 2025}


\maketitle

\begin{abstract}
We use simultaneous Pad\'e approximations to $_3F_2$ hypergeometric functions to estimate from below linear forms in $1$, $\pi\sqrt d$, $\Omega_D/\pi$ and $\pi/\Omega_D$ with integral coefficients, for certain choices of positive integer $d$ and negative integer $D$, where $\Omega_D$ is (the square of) a Chowla--Selberg period attached to the imaginary quadratic field $Q(\sqrt{D})$.
\end{abstract}


This note is not an independent research investigation but rather an addendum to our paper \cite{Zu05} published some twenty years ago.
There we used Ramanujan's and Ramanujan-type hypergeometric formulas for $1/\pi$ and explicit simultaneous Pad\'e approximations to the corresponding hypergeometric functions to produce reasonable estimates for the irrationality measures of $\pi\sqrt d$, for the preselected list $d\in\{1,2,3,10005\}$.
What we have realised now is that we do a bit more in \cite{Zu05}, as some other interesting periods besides the multiples of $\pi$ are approximated as well.
The text below gives some details on the extra.

We consider an integer $D<0$ such that either $D\equiv1\pmod4$, in which case we assume that $D$ is square-free, or $D\equiv0\pmod4$, in which case we assume that $D/4$ is square-free.
Denote by $h=h_D$ the class number of the inaginary quadratic field $\mathbb Q(\sqrt D)$.
Then we define
\[
\Omega_D=\frac{2\pi}{|D|}\Bigg(\prod_{j=1}^{|D|}\Gamma\bigg(\frac j{|D|}\bigg)^{\left(\frac{D}{j}\right)}\Bigg)^{1/h_D},
\]
where $\big(\frac jk\big)$ stands for the Kronecker--Legendre symbol.
Observe that $\Omega_D$ is, up to a positive algebraic multiple, the square of the period constructed in the original paper \cite{SC67} of Selberg and Chowla, but for the purpose of this note it will be convenient to refer to $\Omega_D$ as to the Chowla--Selberg period.

In \cite{Zu05}, we construct Pad\'e type II approximations to $1$ and the functions
\[
f_i(z)=\bigg(z\frac{\d}{\d z}\bigg)^i{}_3F_2\bigg(\begin{matrix} s, \, \frac12, \, 1-s \\ 1, \, 1 \end{matrix}\biggm|z\bigg)
=\sum_{\nu=0}^\infty\frac{(s)_\nu(\frac12)_\nu(1-s)_\nu}{n!^3}\,n^iz^n
\quad\text{for}\; i=0,1,2,
\]
where $s\in\{\frac13,\frac14,\frac16\}$, and apply them at special values of $z$ exclusively for particular $\mathbb Q$-linear combinations of $f_0(z)$ and $f_1(z)$.
These linear combinations are known to evaluate to a number $\sqrt d/\pi$ for some positive integer $d$, thus giving us good rational approximations to the latter.
At the same time we deliberately ignore the fact that our Pad\'e approximations lead to the $\mathbb Q$-linear independence, in a quantitative form, of the \emph{four} numbers $1,f_0(1/Z),f_1(1/Z),f_2(1/Z)$ for integral $Z$ with $|Z|$ sufficiently large.
It happens that for the special values $z=1/Z$ treated in \cite{Zu05} the $\mathbb Q$-linear span of these four numbers coincides with the $\mathbb Q$-linear span of
\begin{equation}
1, \, \frac{\sqrt d\,\Omega_D}{\pi^2}, \, \frac{\sqrt d}{\pi}, \, \frac{\sqrt d}{\Omega_D}
\label{span}
\end{equation}
for a related choice of $D<0$ and $d>0$.
This leads (after rescaling the latter span by $\pi/\sqrt{d}$) to the following theorem, which summarises the actual outcomes of the construction and analysis in~\cite{Zu05}.

\begin{theorem}
\label{main}
For $(D,d)\in\{(-148,1),(-232,2),(-267,3),(-163,10005)\}$, define
\[
\mu_{D,d}=57.53011083\dots, \,
13.93477619\dots, \,
44.12528464\dots, \,
10.02136339\dots,
\]
respectively\textup; these are the estimates for the corresponding irrationality measures of $\pi\sqrt d$ computed in \textup{\cite[Theorem]{Zu05}}.\footnote%
{The class numbers of the quadratic fields $\mathbb Q(\sqrt{-148}),\mathbb Q(\sqrt{-232}),\mathbb Q(\sqrt{-267})$ are equal to~2, while the class number of $\mathbb Q(\sqrt{-163})$ is~1.}
Then for any $\eps>0$ and any collection of integers $m_0,m_1,m_2,m_3$ with $m=\max\{|m_0|,|m_1|,|m_2|,|m_3|\}\ge m^*(\eps)$ we have the following estimate\textup:
\[
|m_0+m_1\pi\sqrt d+m_2\,\Omega_D/\pi+m_3\pi/\Omega_D|\ge m^{1-\mu_{D,d}-\eps}.
\]
\end{theorem}

In the appendix we give the explicit connection between the transcendental entries in \eqref{span} and the values of hypergeometric functions; related algorithms can be found in \cite{Gu21a,Gu21b,Mi22}.
The origin of those expressions rests on the modular parameterisation of the hypergeometric function $f_0(z)$, the closedness of differentiation in the ring of quasimodular forms and the Selberg--Chowla formula \cite{SC67} evaluating the CM values of modular forms as quotiens of the gamma values.
For example, in the `classical' case $s=\frac16$ we have
\[
{}_3F_2\biggl(\begin{matrix} \frac12, \, \frac16, \, \frac56 \\ 1, \, 1 \end{matrix} \biggm|\frac{1728}{j(\tau)}\biggr)
=E_4(\tau)^{1/2},
\]
where $j(\tau)=1728E_4(\tau)^3/(E_4(\tau)^3-E_6(\tau)^2)$ is the modular invariant and
\begin{gather*}
E_2(\tau)=1-24\sum_{m=1}^\infty\frac{mq^m}{1-q^m}, \\
E_4(\tau)=1+240\sum_{m=1}^\infty\frac{m^3q^m}{1-q^m}
\;\;\text{and}\;\;
E_6(\tau)=1-504\sum_{m=1}^\infty\frac{m^5q^m}{1-q^m},
\quad q=e^{2\pi i\tau},
\end{gather*}
are the Eisenstein series. Then $j\bigl((1+\sqrt{-163})/2\bigr)=-1728\cdot53360^3$ and the Selberg--Chowla formula leads to explicit evaluations of $E_2(\tau),E_4(\tau),E_6(\tau)$ at $\tau=(1+\sqrt{-163})/2$ in terms of $\pi$ and~$\Omega_D$. At the same time, $f_0(1/Z),f_1(1/Z),f_2(1/Z)$ for $Z=-53360^3$ received their expressions in terms of the values of the Eisenstein series.
Very similar parameterisations and the same scheme work for $s=\frac13,\frac14$ (and $s=\frac12$ as well, though this case is neither discussed in \cite{Zu05} nor here).

Honestly speaking, the proof of the arithmetic aspects in \cite[Section 2]{Zu05} should be slightly adapted to accommodate the presence of $f_2(z)$ in approximations (originally, only linear combinations of $f_0(z)$ and $f_1(z)$ were taken into account). This modifies the required denominators $D_n$ (in the notation of \cite{Zu05}) but does not affect their asymptotics as given in \cite[Lemmas~3--5]{Zu05}.

More quadruples can be treated for which linear independence estimates can be similarly produced from the general results in \cite{Zu05}; those are linked to the known formulas for $\sqrt d/\pi$ carefully listed in \cite{CG21}.
Furthermore, Theorem~\ref{main} implies estimates for the nonquadraticity measure of $\Omega_D/\pi$ (by setting $m_1=0$).
Finally, notice that very fine bounds for the irrationality measures of the quantities $\omega_D$ related to our Chowla--Selberg periods via the formula $|D|\Omega_D/(2\pi)=\omega_D^2$ are discussed in our recent joint paper \cite{CZ25} with Cohen.

Note that our theorem only discusses the estimates; the fact that the numbers $\Omega_D$ and $\pi$ are algebraically independent is already a consequence of the results of Chudnovsky \cite{Ch76,Wa77}, while Nesterenko's theorem in \cite{Ne96} implies that the three numbers $\Omega_D$, $\pi$ and $\exp(\pi\sqrt{|D|})$ are algebraically independent over the rationals, for any choice of $D<0$ with $D\equiv0,1\pmod4$.

\section*{Appendix}
For $D=-148$, $d=1$ we take $s=\frac14$, $Z=-882^2$. Then
\begin{gather*}
f_0(1/Z)=\frac{42\sqrt d\,\Omega_D}{\pi^2}, \quad
1123f_0(1/Z)+21460f_1(1/Z)=\frac{3528\sqrt d}{\pi}, \\
157655f_0(1/Z)+6024969f_1(1/Z)+115132900f_2(1/Z)=\frac{37044\sqrt d}{\Omega_D}.
\end{gather*}

For $D=-232$, $d=2$ we take $s=\frac14$, $Z=99^4$. Then
\begin{gather*}
f_0(1/Z)=\frac{99\sqrt d\,\Omega_D}{\pi^2}, \quad
4412f_0(1/Z)+105560f_1(1/Z)=\frac{9801\sqrt d}{\pi}, \\
77862889f_0(1/Z)+3725845296f_1(1/Z)+89143308800f_2(1/Z)=\frac{3881196\sqrt d}{\Omega_D}.
\end{gather*}

For $D=-267$, $d=3$ we take $s=\frac13$, $Z=-500^2$. Then
\begin{gather*}
f_0(1/Z)=\frac{75\sqrt d\,\Omega_D}{\pi^2}, \quad
827f_0(1/Z)+14151f_1(1/Z)=\frac{1500\sqrt d}{\pi}, \\
684107f_0(1/Z)+23406555f_1(1/Z)+400501602f_2(1/Z)=\frac{30000\sqrt d}{\Omega_D}.
\end{gather*}

For $D=-163$, $d=10005$ we take $s=\frac16$, $Z=-53360^3$. Then
\begin{gather*}
f_0(1/Z)=\frac{2\sqrt d\,\Omega_D}{\pi^2}, \quad
13591409f_0(1/Z)+545140134f_1(1/Z)=\frac{426880\sqrt d}{\pi}, \\[6pt]
\begin{aligned}
277089597908329f_0(1/Z)+22227667570529352f_1(1/Z) \quad&
\\[-9pt]
+891533297092613868f_2(1/Z)
&
=\frac{136669900800\sqrt d}{\Omega_D}.
\end{aligned}
\end{gather*}

\begin{acknowledgements}
An inspiration for this write-up came from recent discussions with several colleagues on two disjoint topics: modular forms, their periods and CM evaluations\,---\,on one hand, and simultaneous Pad\'e (or Pad\'e type II) approximations\,---\,on the other.
I thank Dmitry Badziahin, Florian Breuer, Francis Brown, Heng Huat Chan, Vesselin Dimitrov, Tiago Jardim da Fonseca, Vasily Golyshev and Pengcheng Zhang for related conversations on these themes.
I further acknowledge the hospitality of the University of Sydney and of the University of Newcastle during my stays in these institutions, in the unusually cold and wet August 2025.

I am grateful for the financial support and hospitality of the Sydney Mathematical Research Institute (SMRI) during July--August 2025.
This work was supported in part from the NWO grant OCENW.M.24.112.
\end{acknowledgements}


\end{document}